\documentclass[12pt]{article}
\usepackage{amssymb,xspace}

\def\header{A.Miller\hfill Recursion Theorem \hfill}
\def\la{\langle}
\def\om{\omega}
\def\proof{\par\noindent Proof\par\noindent}
\def\qed{\par\noindent QED\par\bigskip}
\def\ra{\rangle}
\def\rmand{\mbox{ and }}
\def\self{self-constructing\xspace}
\def\smn{s-m-n\xspace}
\def\st{\;:\;} 
\def\su{\subseteq}

\newtheorem{theorem}{Theorem}
\newtheorem{lemma}{Lemma}

\begin{document}

\begin{center}
{\large The Recursion Theorem and Infinite Sequences}
\end{center}

\begin{flushright}
Arnold W. Miller
\end{flushright}

\def\address{\begin{flushleft}
Arnold W. Miller \\
miller@math.wisc.edu \\
http://www.math.wisc.edu/$\sim$miller\\
University of Wisconsin-Madison \\
Department of Mathematics, Van Vleck Hall \\
480 Lincoln Drive \\
Madison, Wisconsin 53706-1388 \\
\end{flushleft}}

\begin{center}  Abstract  \end{center}
\begin{quote}
In this paper we use the Recursion Theorem to show the
existence of various infinite sequences and sets.  Our main result
is that there is an increasing sequence $e_0<e_1<\cdots$
such that $W_{e_n}=\{e_{n+1}\}$ for every $n$.  
Similarly, we prove that there exists an increasing sequence
such that $W_{e_n}=\{e_{n+1},e_{n+2},\ldots\}$ for every $n$.
We call a nonempty computably enumerable
set $A$ \self if $W_e=A$ for every $e\in A$.  We show that
every nonempty computable enumerable set which is disjoint from
an infinite computable set is one-one equivalent to a \self set.
\end{quote}

Kleene's Recursion Theorem  says the following: 
\begin{quote}
For any computable function
$f$ there exists an $e$ with $\psi_e=\psi_{f(e)}$. 
\end{quote}
In this Theorem
$\la\psi_e:e\in\om\ra$ is a standard computable numbering
of all partial computable functions.  For example, $\psi_e$ might be the
partial function computed by the $e^{th}$ Turing machine.  
The number $e$ is referred to as a fixed point for $f$ and
this theorem is also called the
Fixed Point Theorem.

For a proof of Kleene's Theorem see any of the standard references, Cooper
\cite{cooper}, Odifreddi \cite{odifreddi}, Rogers \cite{rogers}, or Soare
\cite{soare}.  See especially Smullyan \cite{smullyan} for many
variants and generalizations of the fixed
point theorem.   The Recursion Theorem applies to all acceptable
numberings (in the sense of Rogers, see Odifreddi \cite{odifreddi} p.215-221).
All natural enumerations
are acceptable.  We use $W_e$ to denote the domain of $\psi_e$ and
hence $\la W_e:e\in\om\ra$ is a uniform computable listing of
all computably enumerable sets.

The proof of the Recursion Theorem is short but tricky.  It can be uniformized
to yield what is called the Recursion Theorem with Parameters.  The proof also
yields an infinite computable set of fixed points by using the
Padding Lemma. See Soare \cite{soare} pages 36-37.

We will use the following version of the Recursion Theorem with Parameters
which includes a uniform use of the Padding Lemma:

\begin{lemma}\label{lem1}
For any computable function $f:\om\times\om\to\om$ there is
a computable function $h:\om\to\om$ such that 
$W_{h(x)}$ is an infinite computable set for every $x$ and
for every $y\in W_{h(x)}$ we have that
$\psi_{f(x,y)}=\psi_{y}$.
\end{lemma}
The proof is left to the reader.

\bigskip

It is an exercise (see Miller \cite{miller} section 13)
to show (using Smullyan's double recursion theorem
or more particularly its $n$-ary generalization, see \cite{smullyan} Chapter IX)
that for any $n>0$ there is a sequence
 $$e_0<e_1<\cdots<e_n$$
such that $W_{e_i}=\{e_{i+1}\}$ for $i<n$ and 
$W_{e_n}=\{e_0\}$.   It occurred to us to ask if it would be possible
to have an infinite sequence like this.  We show that it is.

Here is our main result:

\begin{theorem}\label{main}
There is a strictly increasing sequence: $$e_0<e_1<\cdots< e_n<\cdots$$ 
such that
$$W_{e_n}=\{e_{n+1}\} \mbox{ for every } n.$$
\end{theorem}
\proof

We use $W_{e,s}$ to denote the set of all $y<s$ such that
$\psi_e(y)$ converges in less than $s$-steps.
We use $\la x,y\ra$ to denote a pairing function,
a computable bijection from $\om^2$ to $\om$, e.g.,
$\la x,y\ra= 2^x(2y+1)-1$.

Let $q(e,x)$ be a computable function such that for all $x$ and $e$:
$$W_{q(e,x)}=\left\{
\begin{array}{ll}
\{y\} & \mbox{ if } (\exists s \exists y\in W_{e,s} \;\; y>x) \mbox{ and 
$\la s,y\ra$ is the least such pair} \\
\emptyset & \mbox{ otherwise. }
\end{array}\right.$$
Such a $q$ is constructed by a standard argument
using the \smn or Parameterization Theorem.
To see this one defines a partial computable function $\theta$ 
as follows:
$$\theta(e,x,y)=\left\{
\begin{array}{ll}
0 & \mbox{ if } (\exists s \exists y\in W_{e,s} \;\; y>x) \mbox{ and 
$\la s,y\ra$ is the least such pair} \\
\uparrow & \mbox{ otherwise. }
\end{array}\right.$$
The uparrow stands for a computation that diverges, i.e, does not halt.
By the \smn Theorem there is a computable $q$ such that
$$\psi_{q(e,x)}(y)=\theta(e,x,y)\mbox{ for all } e,x,y.$$

Using Lemma \ref{lem1} let $h$ be a computable function 
such that for every $e$,
the set $W_{h(e)}$ is an infinite set of fixed points for
$q(e,\cdot)$, i.e., $W_x=W_{q(e,x)}$ for all $x\in W_{h(e)}$.
Let $e$ be a fixed point for $h$, so $W_e=W_{h(e)}$.

Let $x$ be any element of $W_e$.  Then $x\in W_{h(e)}$ so
$W_x=W_{g(x,e)}=\{y\}$ where $y>x$ and $y\in W_e$.  Hence, starting
with any $e_0\in W_e$ we get an infinite increasing sequence as
required.
\qed

Note that to obtain a sequence with $W_{e_{n+1}}=\{e_n\}$ is trivial.
Also note that the sequence in Theorem \ref{main} must be
computable.  This is not necessarily true
for our next result:

\begin{theorem} \label{two}
There exists a computable strictly increasing sequence
$\la e_n:n<\om\ra$ such that for every $n$
$$W_{e_n}=\{e_m\st m>n\}.$$
\end{theorem}
\proof
Using the \smn Theorem find $q$ a computable function
such that for every $x$ and $e$:
$$W_{q(e,x)}=\{\max(W_{e,s})\st s\in \om\}\setminus\{0,1,\ldots,x\}.$$

As in the above proof, let $h$ be a computable function such that for every $e$,
the set $W_{h(e)}$ is an infinite set of fixed points for
$q(e,\cdot)$, i.e., $W_x=W_{q(e,x)}$ for all $x\in W_{h(e)}$.
Let $e$ be a fixed point for $h$, so $W_e=W_{h(e)}$.

Note that $W_e$ is infinite and let 
$$\{e_0<e_1<e_2<\ldots\}=\{\max(W_{e,s})\st s\in \om\}.$$
For any $x\in W_e=W_{h(e)}$
we have that 
$$W_{x}=W_{q(e,x)}=\{\max(W_{e,s})\st s\in \om\}\setminus\{0,1,\ldots,x\}.$$
Hence for any $n$ we have
that
$$W_{e_n}=\{e_m\st m>n\}.$$
\qed

A variation on this theorem would be to get a computable 
strictly increasing sequence
$e_0<e_1<\cdots$ such that 
$$\psi_{e_n}(m)=e_{n+m+1} \mbox{ for every } n,m<\om.$$
The proof of this is left as an exercise for the reader.

\bigskip

Usually the first example given of an application
of the Recursion Theorem is to prove that there exists an
$e$ such that $W_e=\{e\}$.
We say that a nonempty computably enumerable set
$A$ is \self iff for all $e\in A$ we have that
$W_e=A$.  So $W_e=\{e\}$ is an example of a \self set.
Our next result shows there are many \self sets.

\begin{theorem} \label{self}
For any nonempty computably enumerable set $B$ the following are equivalent:
\begin{enumerate}
\item $B$ is disjoint from an infinite computable set.
\item There is a computable permutation $\pi$ of $\om$ such
that $A=\pi(B)$ and $A$ is \self.
\end{enumerate}
\end{theorem}
\proof

$(2)\to (1)$
 
Let $E$ be an infinite computable set such that for
every $e\in E$ we have that $W_e=\emptyset$. Any \self
set is disjoint from $E$.

\bigskip
$(1)\to (2)$

Given any $e$ consider the following computably enumerable set
$Q_e$.  Let
$$\{c_0<c_1<c_2<\ldots\}=\{\max(W_{e,s})\st s\in \om\}$$
which may be finite or even empty.  Then put
$$Q_e=\{c_n \st n\in B\}.$$
By the \smn Theorem we can find a computable $q$ such that
for every $e$:
$$W_{q(e)}=Q_e.$$
By the Padding Lemma there is a computable $h$ such that for
every $e$ the set $W_{h(e)}$ is infinite and for all
$x\in W_{h(e)}$ we have that $W_x=W_{q(e)}$.
Now let $e$ be a fixed point for $h$ so that
$W_{h(e)}=W_e$.   Let $A=Q_e$.   Then for all
$x\in A$ we have that $W_x=A$.  So $A$ is \self.  

To get $\pi$ let $D$ and $E$ be two infinite pairwise 
disjoint computable sets
disjoint from $B$.
Take one-one computable enumerations of them:
$$D=\{d_n\st n<\om\} \rmand E=\{e_n\st n<\om\}.$$
Note that $W_e=W_{h(e)}$ is infinite and let $C$ be the 
infinite computable set:
$$C=\{c_0<c_1<c_2<\ldots\}=\{\max(W_{e,s})\st s\in \om\}.$$
Since $W_e=W_{h(e)}$ is a set of
fixed points, it is coinfinite and hence $C\su W_e$ is
coinfinite.
Take a one-one computable enumeration of the complement $\overline{C}$ of
$C$:
$$\overline{C}=\{\overline{c}_n\st n<\om\}.$$

Now we can define $\pi$:
$$
\pi(c_n)=\left\{
\begin{array}{ll}
d_{2n}   & \mbox{ if } n\in D \\
d_{2n+1} & \mbox{ if } n\in E \\
n        & \mbox{ otherwise }
\end{array}\right.$$ 
$$ \pi(\overline{c}_n)=e_n.\hspace{1.4in}$$
Note that $\pi$ bijectively maps $C$ to $\overline{E}$ and
$\overline{C}$ to $E$.  Furthermore if $n\in B$ then
$\pi(c_n)=n$, and since $A=\{c_n\st n\in B\}$ we
have that $\pi(A)=B$.

\qed

As a corollary we get that there are \self sets of each finite
cardinality and there is a \self set which is not computable,
in fact, there is a creative \self set.

It is not hard to show that
$$S=\{e\st W_e \mbox{ is \self }\}$$
is $\Pi^0_2$-complete.  To see this first note that it is
easy to show that $S$ is $\Pi^0_2$.   We can get
a many-one reduction $f$ of 
$$\mbox{Tot}=^{\rm def}\{e\st W_e=\om\}$$ to $S$ as follows. 
Fix an infinite \self set $A$ with
one-to-one computable enumeration $A=\{a_n\st n\in\om\}$.  By the \smn
Theorem construct a computable $f$ so that for every $e$:
$$W_{f(e)}=\{a_n\st n\in W_e\}.$$
Hence, $e\in$ Tot iff $f(e)\in S$.
Since Tot is $\Pi^0_2$-complete
(see Soare \cite{soare} page 66), $S$ is too.

\address


\newpage

\def\header{ Appendix \hfill Electronic on-line version only\hfill }
\markboth\header\header

\setcounter{page}{1}
\setcounter{theorem}{0}

\begin{center}
 Appendix
\end{center}

\bigskip

The appendix is not intended for final publication but for
the on-line electronic version only.

\bigskip

\begin{theorem}
There exists an infinite strictly increasing sequence
$\la e_n:n<\om\ra$ such that for every $n$
$$W_{e_n}=\{e_m\st m>n\}$$
and the sequence is not computable.
\end{theorem}
\proof
This is a combination of the proof of Theorem \ref{two} and
Theorem \ref{self}.
Fix $B$ a computably enumerable but not computable set.
For any $e$ let
$$\{c_0<c_1<c_2<\ldots\}=\{\max(W_{e,s})\st s\in \om\}$$
By the \smn Theorem find computable $q$ such that for any $x$ and $e$
$$W_{q(e,x)}=\{c_n \st n\in B \rmand c_n>x\}.$$
Take $h$ computable so that $W_{h(e)}$ is an infinite set of
fixed points for $q(e,\cdot)$, i.e., $x\in W_{h(e)}$ implies
$W_x=W_{q(e,x)}$.  Take $e$ to be a fixed point for $h$, so
$W_{h(e)}=W_e$.
Then putting
$$\{e_n \st n<\om\}=\{c_n\st n\in B\}$$
is not computable and have the required property.
\qed

\begin{theorem}
There is a computable
strictly increasing sequence
$e_0<e_1<\cdots$ such that 
$$\psi_{e_n}(m)=e_{n+m+1} \mbox{ for every } n,m<\om.$$
\end{theorem}
\proof
As in the proof of Theorem \ref{two} given
any $e$ let
$$\{e_0<e_1<e_2<\ldots\}=\{\max(W_{e,s})\st s\in \om\}$$
and using the \smn Theorem find $q$ a computable function
such that for every $x$, $m$, and $e$:
$$\psi_{q(e,x)}(m)=
e_{n+m+1} \mbox{ where $n$ is minimal so that $e_n\geq x$}\}.$$

Let $h$ be a computable function such that for every $e$,
the set $W_{h(e)}$ is an infinite set of fixed points for
$q(e,\cdot)$, i.e., $\psi_x=\psi_{q(e,x)}$ for all $x\in W_{h(e)}$.
Let $e$ be a fixed point for $h$, so $W_e=W_{h(e)}$.

Now for any $n$ we have that $e_n\in W_e=W_{h(e)}$ and
so $\psi_{e_n}=\psi_{q(e,e_n)}$ and hence for every $m$:
$$\psi_{e_n}(m)=e_{n+m+1}.$$
\qed


\begin{thebibliography}{99}

\bibitem{cooper} Cooper, S. Barry {\bf Computability theory}. Chapman \&
Hall/CRC, Boca Raton, FL, 2004. x+409 pp. ISBN: 1-58488-237-9 

\bibitem{miller} Miller, Arnold W. {\bf Lecture Notes in Computability Theory}.
eprint Jan 08:
http://www.math.wisc.edu/$\sim$miller/res

\bibitem{odifreddi} Odifreddi, Piergiorgio { \bf Classical recursion theory. The
theory of functions and sets of natural numbers.} With a foreword by G. E.
Sacks. Studies in Logic and the Foundations of Mathematics, 125. North-Holland
Publishing Co., Amsterdam, 1989. xviii+668 pp. ISBN: 0-444-87295-7 

\bibitem{rogers} Rogers, Hartley, Jr. {\bf Theory of recursive functions and
effective computability.} McGraw-Hill Book Co., New York-Toronto, Ont.-London
1967 xx+482 pp. 

\bibitem{smullyan} Smullyan, Raymond M. {\bf Recursion theory for
metamathematics.} Oxford Logic Guides, 22. The Clarendon Press, Oxford University
Press, New York, 1993. xvi+163 pp. ISBN: 0-19-508232-X 


\bibitem{soare}
Soare, Robert I. {\bf Recursively enumerable sets and
degrees. A study of computable functions and computably generated sets.}
Perspectives in Mathematical Logic. Springer-Verlag, Berlin, 1987. xviii+437 pp.
ISBN: 3-540-15299-7 


\end{thebibliography}
\end{document}